\documentclass[12pt,twoside]{article}
\usepackage{amsmath}

\date{}
\begin{document}

%\centerline{\bf Bulletin of the Malaysian Mathematical Sciences
%Society}

%\centerline{Vol. x, 200x, no. xxx - xxx}

\centerline{}

\centerline{\Large{\bf Operations on fuzzy ideals of
$\Gamma-$semirings}}

\centerline{}

\centerline{\bf{T.K.Dutta}}

\centerline{Department of Pure Mathematics,}

\centerline{Calcutta University, Kolkata, India}

\centerline{E-mail: duttatapankumar@yahoo.co.in}

\centerline{}

\centerline{\bf{Sujit Kumar Sardar}}

\centerline{Department of Mathematics,}

\centerline{Jadavpur University, Kolkata, India}

\centerline{E-mail: (i) sksardar@math.jdvu.ac.in, (ii)
sksardarjumath@gmail.com}

\centerline{}

\centerline{\bf{Sarbani Goswami}}

\centerline{Lady Brabourne College, Kolkata, India}

\centerline{E-mail: sarbani7$_{-}$goswami@yahoo.co.in}

\newtheorem{Theorem}{\quad Theorem}[section]

\newtheorem{Definition}[Theorem]{\quad Definition}

\newtheorem{Corollary}[Theorem]{\quad Corollary}

\newtheorem{Lemma}[Theorem]{\quad Lemma}

\newtheorem{Example}[Theorem]{\emph{Example}}

\newtheorem{Proposition}[Theorem]{Proposition}

\begin{abstract}
 The purpose of this paper is to introduce different types of operations on fuzzy ideals of
 $\Gamma-$semirings and to prove subsequently that these oprations
 give rise to different structures such as complete lattice, modular
 lattice on some restricted class of fuzzy ideals of
 $\Gamma-$semirings. A characterization of a regular
 $\Gamma-$semiring has also been obtained in terms of fuzzy subsets.
\end{abstract}


\section{Introduction} If we remove the restriction of having
additive inverse of each element in a ring then a new algebraic
structure is obtained what we call a semiring. Semiring has found
many applications in various fields. In this regard we may refer to
Golan's \cite{re:Golan} and Weinert's \cite{re:Hebisch} monographs.
Semiring arises very naturally as the nonnegative cone of a totally
ordered ring. But the nonpositive cone of a totally ordered ring
fails to be a semiring because the multiplication is no longer
defined. One can provide an algebraic home, called
$\Gamma-$semiring, to the nonpositive cone of a totally ordered
ring. The notion of $\Gamma-$semiring was introduced by M.M.K.Rao
\cite{re:Rao} in 1995 as a generalization of semiring as well as of
$\Gamma-$ring. Subsequently by introducing the notion of operator
semirings of a $\Gamma-$semiring Dutta and Sardar enriched the
theory of $\Gamma-$semirings. In this connection we may refer to
\cite{re:Dutta}. The motivation for this paper is the fact that
$\Gamma-$semiring is a generalization of semiring as well as of
$\Gamma-$ring and fuzzy concepts of Zadeh \cite{re:Zadeh} has been
successfully applied to $\Gamma-$rings and semirings by Jun et al
\cite{re:Jun} and Dutta et al \cite{re:Chanda}, \cite{re:Biswas}. We
define  here some compositions of fuzzy ideals in a
$\Gamma-$semiring and study the structures of the set of fuzzy
ideals of a $\Gamma-$semiring. Among other results we have deduced
that sets of fuzzy left ideals and fuzzy right ideals form a
zero-sum free semiring with infinite element. We have also deduced
that fuzzy ideals of a $\Gamma-$semiring is a complete lattice which
is modular if every fuzzy ideal is a fuzzy k-ideal.

\section{Preliminaries}

\begin{Definition} \textnormal{\cite{re:Rao}}Let S and $\Gamma$ be two additive
commutative semigroups. Then S is called a $\Gamma-$semiring if
there exists a mapping $\\ S \times \Gamma \times S \rightarrow S $
(images to be denoted by $a \alpha b $ for $a, b \in S $ and $\alpha
\in \Gamma) $ satisfying the following conditions:\\
(i) $(a+b) \alpha c =a \alpha c+ b \alpha c$,\\
(ii) $ a \alpha (b+c)= a \alpha b+a \alpha c $,\\
(iii) $a (\alpha+\beta)b=a \alpha b +a \beta b$,\\
(iv) $ a \alpha (b \beta c)=(a \alpha b) \beta c $
for all $a,b,c \in S $ and for all $\alpha, \beta \in \Gamma $.\\
\\Further, if in a $\Gamma-$semiring, $(S,+)$ and $(\Gamma,+)$ are both
monoids and\\ (i) $ 0_{S} \alpha x=0_{S} =x \alpha 0_{S}$\\(ii) $ x
0_{\Gamma} y=0_{S}=y 0_{\Gamma} x$ for all $x,y \in S$ and for all
$\alpha \in \Gamma $ then we say that S is a $\Gamma-$semiring with
zero.\end{Definition}
 Throughout this paper we consider $\Gamma-$semiring with zero. For
simplification we write 0 instead of $0_{S}$ and $0_{\Gamma}$ which
will be clear from the context.

 \begin{Definition} \textnormal{\cite{re:Zadeh}} Let S be a non empty set. A mapping
 $\mu: S \rightarrow [0,1]$ is called a fuzzy subset of S.\end{Definition}

\begin{Definition} \textnormal{\cite{re:Goswami}}Let $\mu$ be a non empty fuzzy subset of a
$\Gamma-$semiring S (i.e. $\mu(x) \neq 0 $ for some $x \in S
 $). Then $\mu $ is called a fuzzy left ideal [ fuzzy right ideal]
 of S if \\(i) $\mu(x+y) \geq min [\mu(x), \mu(y)]$ and
  \\~~~~~~~~(ii) $\mu(x \gamma y) \geq \mu(y) $ [resp. $\mu(x \gamma y) \geq
  \mu(x)]$ for all $ x,y \in S, \gamma \in \Gamma$.\end{Definition}

A fuzzy ideal of a $\Gamma-$semiring S is a non empty fuzzy subset
of S which is a fuzzy left ideal as well as a fuzzy right ideal of
S.

\begin{Definition} \textnormal{\cite{re:Golan}}Let S be a non empty set and `+' and `.' be
two binary operations on S, called addition and multiplication
respectively. Then $( S, + , .)$ is called a hemiring (resp.
semiring) if

(i) (S, +) is a commutative monoid with identity element 0;

(ii) (S, .) is a semigroup (resp. monoid with identity element 1);

(iii) $ a.(b+c)=a.b+a.c $ and $ (b+c).a=b.a+c.a $ for all $ a, b, c
\in S $.

(iv) $ a.0=0.a=0 $ for all $ a \in S $;

(v) $1 \neq 0 $.

A hemiring S is said to be zero-sum free if $ a+b=0$ implies that
$a=b=0  $ for all $a,b \in S$.

An element $a$ of a hemiring S is infinite iff $a+s=a$ for all $s
\in S$.
\end{Definition}
For more on preliminaries we may refer to the references and their
references.

\section{Operations on fuzzy ideals}

 Throughout this paper unless otherwise mentioned S denotes a $\Gamma$-semiring with
 unities\textnormal{\cite{re:Dutta}} and $ FLI(S) $, $ FRI(S) $ and $ FI(S) $
denote respectively the set of all fuzzy left ideals, the set of all
fuzzy right ideals and the set of all fuzzy ideals of the $ \Gamma
$-semiring S. Also in this section we assume that $\mu(0)=1$ for a
fuzzy left ideal (fuzzy right ideal, fuzzy ideal) $\mu$ of a
$\Gamma-$semiring ($\Gamma-$hemiring) S.

\begin{Definition} Let S be a $ \Gamma $-semiring and $
\mu_{1}, \mu_{2} \in FLI(S) ~[ FRI(S),~FI(S)] $. Then the sum $
\mu_{1}\oplus \mu_{2} $, product $ \mu_{1}\Gamma \mu_{2} $ and
composition $ \mu_{1}\circ \mu_{2} $ of $ \mu_{1} $ and $\mu_{2} $
are defined as follows:

$$ (\mu_{1} \oplus \mu_{2})(x)= \sup_{x=u+v}[\min[ \mu_{1}(u),
\mu_{2}(v) ]: u, v \in S]$$
\\ $ ~~~~~~~~~~~~~~~~~~~~~~~~~~~~~~= 0 $ if
for any $ u,v \in S, u+v \neq x $.

$$ ( \mu_{1} \Gamma \mu_{2})(x) = \sup_{x=u \gamma v} [ \min [
\mu_{1}(u), \mu_{2}(v)] :u,v \in S ; \gamma \in \Gamma] $$
\\ $~~~~~~~~~~~~~~~~~~~~~~~= 0$ if for any $ u,v \in S $ and for any $ \gamma
\in \Gamma, ~  u \gamma v \neq x $.

$ ( \mu_{1} \circ \mu_{2})(x)=\displaystyle{ \sup_{ x=
\displaystyle{\sum_{i=1}^{n}} u_{i} \gamma_{i} v_{i}}}
[\displaystyle{ \min_{ 1\leq i \leq n }} [ \min [ \mu_{1}(u_{i}),
\mu_{2}(v_{i})]]:  u_{i}, v_{i} \in S, \gamma_{i} \in \Gamma ] $

~~~~~~~~~~~~~~~~~= 0 otherwise.\end{Definition}

\textbf{ Note.} Since S contains $ 0 $, in the above definition the
case $ x \neq u+v $ for any $ u,v \in S $ does not arise. Similarly
since S contains left and right unity, the case $ x \neq
\displaystyle{\sum_{i}} u_{i} \gamma_{i} v_{i} $ for any $
u_{i},v_{i} \in S, \gamma_{i} \in \Gamma $ does not arise. In case
of product of $\mu_{1}$ and $\mu_{2}$ if S has strong left or right
unity [i.e., there exists $ e \in S, \delta \in \Gamma$ such that $e
\delta a =a\textmd{ for all }a \in S $] then the case $x \neq u
\gamma v $ for any $u,v \in S$ and for any $\gamma \in \Gamma$ does
not arise.
 i.e., in otherwords there are $u,v \in S$ and $\gamma \in \Gamma$ such that $x = u \gamma v $.   \\

\begin{Proposition} Let $ \mu_{1}, \mu_{2} \in FLI(S) [ FRI(S),
FI(S) ] $. Then

$ \mu_{1} \oplus \mu_{2} \in FLI(S) [ $ resp. $ FRI(S), FI(S) ]
$.\label{prop:3.2}\end{Proposition}

\textsl{ Proof.} $(\mu_{1} \oplus
\mu_{2})(0)=\displaystyle{\sup_{0=u+v}[\min[\mu_{1}(u), \mu_{2}(v)]:
u,v \in S]} \\~~~~~~~~~~~~~~~~~~~~~~~\geq
\min[\mu_{1}(0),\mu_{2}(0)]=1 \neq 0 $.
\\Thus $\mu_{1} \oplus
\mu_{2}$ is non empty and $ ( \mu_{1} \oplus \mu_{2})(0)= 1 $.
\\Let $ x, y \in S $ and $ \gamma \in \Gamma $.
Then\\
$ ( \mu_{1} \oplus \mu_{2})(x+y) = \displaystyle{\sup_{x+y= p+q}} [
\min [ \mu_{1}(p), \mu_{2}(q)] : p,q \in S
]\\~~~~~~~~~~~~\geq \displaystyle{\sup_{\begin{array}{l} x= u+v \\
y= s+t
\end{array} }}[ \min [ \mu_{1}(u+s),
\mu_{2}(v+t)]$:  $ u,v,s,t \in S ]\\~~~~~~~~~~~~~~ \geq
\displaystyle{\sup _{\begin{array}{l} x= u+v \\ y= s+t
\end{array} }}[ \min[ \min [ \mu_{1}(u), \mu_{1}(s)], \min [
\mu_{2}(v), \mu_{2}(t)]]$:  $ u,v,s,t \in S ]
\\~~~~~~~~~~~~~ =
\displaystyle{\sup _{\begin{array}{l} x= u+v \\ y= s+t
\end{array} }}[ \min[ \min [ \mu_{1}(u), \mu_{2}(v)], \min [
\mu_{1}(s), \mu_{2}(t)]]$:  $ u,v,s,t \in S ]\\~~~~~~~~~~~~~
=\displaystyle{ \min [ \sup_{x=u+v} [ \min [ \mu_{1}(u),
\mu_{2}(v)]], \sup_{y=s+t} [ \min [ \mu_{1}(s),
\mu_{2}(t)]]]}\\~~~~~~~~~~~~~ = \min [ (\mu_{1} \oplus \mu_{2})(x),
(\mu_{1} \oplus \mu_{2})(y)] $.

Again $ ( \mu_{1} \oplus \mu_{2})(x \gamma y) =\displaystyle{\sup_{
x \gamma y= p+q} [\min [ \mu_{1}(p),
\mu_{2}(q)]]}\\~~~~~~~~~~~~~~~~~~~~~~~~~~~~~~~~\geq
\displaystyle{\sup_{ y=u+v} [\min [ \mu_{1}(x \gamma u), \mu_{2}(x
\gamma v)]]}$
\\$~~~~~~~~~~~~~~~~~~~~~~~~~~~~~~~~~~~~~~~~~~~~~~~~~~~~~~~~$[ Since $ x \gamma y
= x \gamma (u+v) = x \gamma u + x \gamma v ]
\\~~~~~~~~~~~~~~~~~~~~~~~~~~~~~~~~\geq \displaystyle{\sup_{y=u+v} [ \min [
\mu_{1}(u), \mu_{2}(v)]]} = ( \mu_{1} \oplus \mu_{2})(y) $.

Hence $ \mu_{1} \oplus \mu_{2} \in FLI(S) $.

\begin{Proposition} Let $ \mu_{1}, \mu_{2}, \mu_{3} \in
FLI(S) [ FRI(S), FI(S) ] $. Then

(i) $ \mu_{1} \oplus \mu_{2} =  \mu_{2} \oplus \mu_{1}$.

(ii) $ ( \mu_{1} \oplus \mu_{2}) \oplus \mu_{3} = \mu_{1} \oplus (
\mu_{2} \oplus \mu_{3}).$

(iii) $ \theta \oplus \mu_{1} = \mu_{1} = \mu_{1} \oplus \theta $
where $ \theta $ is a fuzzy ideal of S, defined by,

$ \theta(x)= \left \{\begin{array}{l} 1 ~~~\textmd{if}~~x=0 \\ 0
~~~\textmd{if}~~ x \neq 0
\end{array} \right .  $

(iv) $ \mu_{1} \oplus \mu_{1} =\mu_{1} $.

(v) $ \mu_{1} \subseteq \mu_{1} \oplus \mu_{2} $ and

(vi) $ \mu_{1} \subseteq \mu_{2}$ implies that $ \mu_{1} \oplus
\mu_{3} \subseteq \mu_{2} \oplus
\mu_{3}.$\label{prop:3.3}\end{Proposition}

\textsl{ Proof.} (i) We leave it as it follows easily.

(ii) Let $ x \in S $.

$((\mu_{1} \oplus \mu_{2}) \oplus \mu_{3})(x) =\displaystyle{
\sup_{x=u+v}} [ \min [ ( \mu_{1} \oplus \mu_{2})(u), \mu_{3}(v)]: u,
v \in S ]
\\~~~~~~~~~~~~~~~~~= \displaystyle{\sup_{ x=u+v}} [\min [ \displaystyle{\sup_{ u= p+q}} [ \min [ \mu_{1}(p),
\mu_{2}(q)]: p,q \in S ]], \mu_{3}(v)]
\\~~~~~~~~~~~~~~~~=\displaystyle{\sup_{x=u+v}}~~ \displaystyle{\sup_{u=p+q}} [ \min [ \min [
\mu_{1}(p), \mu_{2}(q)], \mu_{3}(v)]] \\~~~~~~~~~~~~~~~~=
\displaystyle{\sup_{ x= p+q+v}} [ \min[ \mu_{1}(p), \mu_{2}(q),
\mu_{3}(v)]] $.

Similarly we can deduce that $ ( \mu_{1} \oplus ( \mu_{2} \oplus
\mu_{3}))(x) =\displaystyle{\sup_{x=p+q+v}}[\min[ \mu_{1}(p),
\mu_{2}(q), \mu_{3}(v)]]$.

Therefore $ ( \mu_{1} \oplus \mu_{2}) \oplus \mu_{3}= \mu_{1} \oplus
( \mu_{2} \oplus \mu_{3}). $

(iii) For any $ x \in S ,\\ ( \theta \oplus \mu_{1})(x)=
\displaystyle{\sup_{x=u+v}}[ \min[ \theta(u), \mu_{1}(v)]$, for $
u,v \in S]\\~~~~~~~~~~~~~ =\min[ \theta(0), \mu_{1}(x)]=
\mu_{1}(x)$.\\Thus $ \theta \oplus \mu_{1}= \mu_{1}$. From (i) $
\mu_{1} \oplus \theta=\theta \oplus \mu_{1}= \mu_{1} $.

(iv) Let $ x \in S $. Then $\\  ( \mu_{1} \oplus \mu_{1})(x) =
\displaystyle{\sup_{x=u+v}} [\min[\mu_{1}(u), \mu_{1}(v)],$ for $ u,
v \in S]
\\~~~~~~~~~~~~~~\leq \displaystyle{\sup_{x=u+v}} \mu_{1}(u+v)=
\mu_{1}(x) $

So $ \mu_{1} \oplus \mu_{1} \subseteq \mu_{1} $

Again $ \mu_{1}(x)=\min[\mu_{1}(0), \mu_{1}(x)]\\~~~~~~~~~~~~~~ \leq
\displaystyle{\sup_{x=u+v}} [\min[\mu_{1}(u), \mu_{1}(v)],$ for $ u,
v \in S] =(\mu_{1} \oplus \mu_{1})(x)$.

Therefore $ \mu_{1} \subseteq \mu_{1} \oplus \mu_{1} $.
Consequently, $ \mu_{1}= \mu_{1} \oplus \mu_{1} $.

(v) Let $ x \in S $. Then $\\  ( \mu_{1} \oplus \mu_{2})(x) =
\displaystyle{\sup_{x=u+v}} [\min[\mu_{1}(u), \mu_{2}(v)],$ for $ u,
v \in S]
\\~~~~~~~~~~~~~~\geq \min[ \mu_{1}(x), \mu_{2}(0)]
= \mu_{1}(x)$.

Thus $ \mu_{1} \subseteq \mu_{1} \oplus \mu_{2} $.

(vi) Let $ \mu_{1} \subseteq \mu_{2}.$ and $ x \in S $. Then $ \\
(\mu_{1} \oplus \mu_{3})(x)=
\displaystyle{\sup_{x=u+v}}[\min[\mu_{1}(u), \mu_{3}(v)],$ for $ u,
v \in S] \\~~~~~~~~~~~~\leq
\displaystyle{\sup_{x=u+v}}[\min[\mu_{2}(u), \mu_{3}(v)],$ for $ u,
v \in S] =(\mu_{2} \oplus \mu_{3})(x).$

Hence $ \mu_{1} \oplus \mu_{3} \subseteq \mu_{2} \oplus \mu_{3}$.

\begin{Proposition} Let $ \mu_{1}, \mu_{2} \in FLI(S) [
FRI(S) , FI(S)] $. Then

$\mu_{1} \circ \mu_{2} \in FLI(S)$ [resp. FRI(S),
FI(S)].\label{prop:3.4}\end{Proposition}

\textsl{Proof.} Since $(\mu_{1} \circ \mu_{2}
)(0)\\~~~~~~~~~~~~~~~~~=\displaystyle{\sup_{0=\displaystyle{\sum^{n}_{i=1}
u_{i} \gamma_{i} v_{i}}}[\min_{1 \leq i \leq n}}
[\min[\mu_{1}(u_{i}), \mu_{2}(v_{i})]] : u_{i}, v_{i} \in S,
\gamma_{i} \in \Gamma, n \in Z^{+}]
\\~~~~~~~~~~~~~~~\geq \min[\mu_{1}(0),
\mu_{2}(0)] =1 \neq 0 $ [Since $\mu_{1}(0)=\mu_{2}(0)=1$],\\ it
follows that $\mu_{1} \circ \mu_{2}$ is nonempty and $(\mu_{1} \circ
\mu_{2} )(0)=1$.
\\Now, for any $x,y \in S , \\ (\mu_{1} \circ \mu_{2}) (x+y)\\~~~= \displaystyle{\sup_{x+y=\displaystyle{\sum^{n}_{i=1} u_{i}
\gamma_{i} v_{i}}}[\min_{1 \leq i \leq n}} [\min[\mu_{1}(u_{i}),
\mu_{2}(v_{i})]] : u_{i}, v_{i} \in S, \gamma_{i} \in \Gamma, n \in
Z^{+}]\\~~~\geq \sup[\min_{\begin{array}{l} 1\leq i \leq m
\\ 1 \leq k \leq l
\end{array} }[\min[\min[\mu_{1}(u_{i}),\mu_{2}(v_{i})],\min[\mu_{1}(p_{k}),\mu_{2}(q_{k})]]]
:\\~~~~~~~~~~~~~x=\displaystyle{\sum^{m}_{i=1} u_{i} \gamma_{i}
v_{i}},y=\displaystyle{\sum^{l}_{k=1} p_{k} \gamma_{k} q_{k}},
u_{i},v_{i},p_{k},q_{k} \in S; \gamma_{i} \in \Gamma; m,l \in
Z^{+}]\\=\min[\displaystyle{\sup_{x=\displaystyle{\sum^{m}_{i=1}
u_{i} \gamma_{i} v_{i}}}[\min_{1 \leq i \leq m}}
[\min[\mu_{1}(u_{i}), \mu_{2}(v_{i})]] : u_{i}, v_{i} \in S,
\gamma_{i} \in \Gamma, m \in Z^{+}],
\\ ~~~~~~~~~~\displaystyle{\sup_{y=\displaystyle{\sum^{l}_{k=1} p_{k} \gamma_{k}
v_{k}}}[\min_{1 \leq k \leq l}} [\min[\mu_{1}(p_{k}),
\mu_{2}(q_{k})]] : p_{k}, q_{k} \in S, \gamma_{k} \in \Gamma, l \in
Z^{+}] ] \\ =\min[(\mu_{1} \circ \mu_{2})(x),(\mu_{1} \circ
\mu_{2})(y)]$.
\\Now $(\mu_{1} \circ \mu_{2})(x \gamma y)\\~~~=\displaystyle{\sup_{x \gamma y=\displaystyle{\sum^{n}_{i=1} u_{i}
\gamma_{i} v_{i}}}[\min_{1 \leq i \leq n}} [\min[\mu_{1}(u_{i}),
\mu_{2}(v_{i})]] : u_{i}, v_{i} \in S, \gamma_{i} \in \Gamma,n \in
Z^{+}]
\\~~~\geq
\displaystyle{\sup_{y=\displaystyle{\sum^{m}_{j=1} s_{j} \delta_{j}
t_{j}}}[\min_{1 \leq j \leq m}} [\min[\mu_{1}(x \gamma s_{j}),
\mu_{2}(t_{j})]] ]\\~~~\geq
\displaystyle{\sup_{y=\displaystyle{\sum^{m}_{j=1} s_{j} \delta_{j}
t_{j}}}[\min_{1 \leq j \leq m}} [\min[\mu_{1}(s_{j}),
\mu_{2}(t_{j})]] ]=(\mu_{1} \circ \mu_{2})(y) $
\\Hence $\mu_{1} \circ \mu_{2} \in FLI(S)$
\begin{Proposition} Let $ \mu_{1}, \mu_{2} \in FLI(S) [
FRI(S), FI(S) ] $. Then $\\ \mu_{1} \Gamma \mu_{2} \subseteq \mu_{1}
\circ \mu_{2} $.\label{prop:3.5}\end{Proposition}

\textsl{ Proof.} If for any $ u,v \in S $ and for any $ \gamma \in
\Gamma, ~  u \gamma v \neq x $ then $\\ \mu_{1} \Gamma \mu_{2}
\subseteq \mu_{1} \circ \mu_{2} $.
\\ Now for any $ x \in S, (\mu_{1}
\circ \mu_{2})(x)=\\~~~~~
\displaystyle{\sup_{x=\displaystyle{\sum^{n}_{i=1} u_{i} \gamma_{i}
v_{i}}}[\min_{1 \leq i \leq n}} [\min[\mu_{1}(u_{i}),
\mu_{2}(v_{i})]] : u_{i}, v_{i} \in S, \gamma_{i} \in \Gamma, n \in
Z^{+}]
\\~~~~~ \geq \displaystyle{\sup_{x=u
\gamma v}}[\min [ \mu_{1}(u), \mu_{2}(v)]]=(\mu_{1} \Gamma
\mu_{2})(x) $. \\Thus $ \mu_{1} \Gamma \mu_{2} \subseteq \mu_{1}
\circ \mu_{2} $.

\begin{Proposition} Let $ \mu_{1}$ be a fuzzy right ideal and
$\mu_{2}$ be a fuzzy left ideal of S. Then $\mu_{1} \Gamma \mu_{2}
\subseteq \mu_{1} \cap \mu_{2} $.\label{prop:3.6}\end{Proposition}

\textsl{ Proof.} Let $ \mu_{1}$ be a fuzzy right ideal and $\mu_{2}$
be a fuzzy left ideal of S. For $x \in S$,
\\$(\mu_{1} \Gamma \mu_{2})(x)=\displaystyle{\sup_{x=u \gamma v}}[\min [
\mu_{1}(u), \mu_{2}(v)]: u,v \in S] \leq \displaystyle{\sup_{x=u
\gamma v}}[\min [ \mu_{1}(u \gamma v), \mu_{2}(u\gamma
v)]]\\~~~~~~~~~~~~~~~\leq \displaystyle{\sup_{x=u \gamma v}}(\mu_{1}
\cap \mu_{2})(u\gamma v)=(\mu_{1} \cap \mu_{2})(x)$.
\\Thus $\mu_{1} \Gamma \mu_{2}
\subseteq \mu_{1} \cap \mu_{2} $.
\\\\The following is a characterization of a regular
 $\Gamma-$semiring in terms of fuzzy subsets.
\begin{Theorem} A $\Gamma-$semiring S is multiplicatively
regular\textnormal{\cite{re:Rao}} if and only if  $\mu_{1} \Gamma
\mu_{2} =\mu_{1} \cap \mu_{2} $ for every fuzzy right ideal
$\mu_{1}$ and every fuzzy left ideal $\mu_{2}$ of
S.\label{prop:3.7}\end{Theorem}

\textsl{Proof.} Let S be a multiplicatively regular
$\Gamma-$semiring and $ \mu_{1}$ be a fuzzy right ideal and
$\mu_{2}$ be a fuzzy left ideal of S. Then by Proposition
\ref{prop:3.6}, \\$\mu_{1} \Gamma \mu_{2} \subseteq \mu_{1} \cap
\mu_{2} $. Let $c \in S$. Since S is multiplicatively regular, there
exists an element $x$ in S and $\gamma_{1}, \gamma_{2} \in \Gamma$
such that $c=c \gamma_{1} x \gamma_{2} c $.
\\Now $(\mu_{1} \Gamma \mu_{2})(c)=\displaystyle{\sup_{c=a
\gamma b}}[min[\mu_{1}(a),\mu_{2}(b)]:a,b \in S; \gamma \in
\Gamma]\\~~~~~~~~~~~~~~~~~~~~~\geq min[\mu_{1}(c \gamma_{1}
x),\mu_{2}(c)]~~~~$ [Since $c=(c \gamma_{1} x) \gamma_{2} c $]
\\$~~~~~~~~~~~~~~~~~~~~~\geq min[\mu_{1}(c),\mu_{2}(c)]=(\mu_{1} \cap
\mu_{2})(c)$.
\\Therefore $(\mu_{1} \cap
\mu_{2}) \subseteq \mu_{1} \Gamma \mu_{2}$ and hence $\mu_{1} \Gamma
\mu_{2} =\mu_{1} \cap \mu_{2} $.
\\Conversely, let S is a $\Gamma-$semiring and for every fuzzy right ideal $\mu_{1}$ and every fuzzy left
ideal $\mu_{2}$ of S, $\mu_{1} \Gamma \mu_{2} =\mu_{1} \cap \mu_{2}
$. Let L and R be a left ideal and a right ideal of S respectively
and let $x \in L \cap R $.
\\So $\lambda_{L}(x)=1=\lambda_{R}(x)$. Thus
$(\lambda_{L}\cap\lambda_{R})(x)=1$. Now since $\lambda_{R} \Gamma
\lambda_{L}=\lambda_{R} \cap \lambda_{L}$, so $(\lambda_{R} \Gamma
\lambda_{L})(x)=1$. Therefore $\displaystyle{\sup_{x=y \gamma z}}
[min[\lambda_{R}(y),\lambda_{L}(z)]:y,z \in S; \gamma \in
\Gamma]=1$.
\\Thus there exists some $r,s \in S$ and $\gamma_{1} \in \Gamma$
such that $\lambda_{L}(s)=1=\lambda_{R}(r)$ for $x=r \gamma_{1} s$.
Then $r \in R$ and $s \in L$ and so $x=r \gamma_{1} s \in R \Gamma
L$. Therefore $L \cap R \subseteq R \Gamma L $. Also $L \cap R
\supseteq R \Gamma L $. Thus $R \Gamma L=R \cap L $. Consequently, S
is multiplicatively regular.

\begin{Proposition} Let $ \mu_{1}, \mu_{2} \in FI(S) $. Then
$$ \mu_{1} \Gamma \mu_{2} \subseteq \mu_{1} \circ \mu_{2} \subseteq
\mu_{1} \cap \mu_{2} \subseteq \mu_{1}, \mu_{2}.
$$\label{prop:3.8}\end{Proposition}

\textsl{Proof.} By Proposition \ref{prop:3.5}, $ \mu_{1} \Gamma
\mu_{2} \subseteq \mu_{1} \circ \mu_{2} $. For any $ x \in S $, if
\\$( \mu_{1} \circ \mu_{2})(x)=0 $ then obviously $\mu_{1} \circ
\mu_{2} \subseteq \mu_{1} \cap \mu_{2}$. Now for any $ x \in S $,
 \\ $( \mu_{1} \circ \mu_{2})(x)\\~~~= \displaystyle{\sup_{x=\displaystyle
{\sum^{n}_{i=1}} u_{i} \gamma_{i} v_{i} }}[\min_{1 \leq i \leq n}
[\min[\mu_{1}(u_{i}), \mu_{2}(v_{i})]] : u_{i}, v_{i} \in S,
\gamma_{i} \in \Gamma, n \in Z^{+}]
\\~~~ \leq \displaystyle{\sup_{x=\displaystyle{\sum^{n}_{i=1}} u_{i} \gamma_{i} v_{i}}}
[\min_{1 \leq i \leq n} [\min[\mu_{1}(u_{i} \gamma_{i} v_{i}),
\mu_{2}(u_{i} \gamma_{i} v_{i})]] : u_{i}, v_{i} \in S, \gamma_{i}
\in \Gamma, n \in Z^{+}] \\~~~\leq \min [ \mu_{1}(x), \mu_{2}(x)] =
(\mu_{1} \cap \mu_{2})(x)$.
\\Therefore $ \mu_{1} \circ \mu_{2} \subseteq \mu_{1} \cap \mu_{2} $.
Again $ (\mu_{1} \cap \mu_{2})(x)=\min [\mu_{1}(x),\mu_{2}(x)] \leq
\mu_{1}(x) $. Thus $ \mu_{1} \cap \mu_{2} \subseteq \mu_{1} $.
Similarly it can be shown that $ \mu_{1} \cap \mu_{2} \subseteq
\mu_{2} $. Hence the proposition.

\begin{Proposition} Let $ \mu_{1}, \mu_{2}, \mu_{3} \in
FLI(S)[ FRI(S), FI(S)] $. Then \\$ \mu_{1} \Gamma \mu_{2} \subseteq
\mu_{3} $ if and only if $ \mu_{1} \circ \mu_{2} \subseteq \mu_{3}
$.\label{prop:3.9}\end{Proposition}

\textsl{Proof.} Since $ \mu_{1} \Gamma \mu_{2} \subseteq \mu_{1}
\circ \mu_{2}$ it follows that $ \mu_{1} \circ \mu_{2} \subseteq
\mu_{3} $ implies that

$ \mu_{1} \Gamma \mu_{2} \subseteq \mu_{3} $. Assume that $ \mu_{1}
\Gamma \mu_{2} \subseteq \mu_{3}$. Let $ x \in S $ and\\ $
x=\displaystyle{\sum_{i=1}^{n}} u_{i} \gamma_{i} v_{i}, u_{i}, v_{i}
\in S, \gamma_{i} \in \Gamma, n \in Z^{+} $.
\\Then $ \mu_{3}(x) = \mu_{3}( \displaystyle{\sum_{i=1}^{n}} u_{i} \gamma_{i} v_{i})
\\~~~~~~~~~~~~~~~~\geq \min[\mu_{3}(u_{1} \gamma_{1} v_{1}), \mu_{3}(u_{2} \gamma_{2} v_{2}),
........., \mu_{3}(u_{n} \gamma_{n} v_{n})] \\~~~~~~~~~~~~~~~~\geq
\min[ (\mu_{1} \Gamma \mu_{2})( u_{1} \gamma_{1} v_{1}),(\mu_{1}
\Gamma \mu_{2})( u_{2} \gamma_{2} v_{2}),.........,(\mu_{1} \Gamma
\mu_{2})(u_{n} \gamma_{n} v_{n})] \\~~~~~~~~~~~~~~~~\geq
\min[\min[\mu_{1}(u_{1}),
\mu_{2}(v_{1})],.......,\min[\mu_{1}(u_{n}), \mu_{2}(v_{n})] $.

$ \mu_{3}(x) \geq
\displaystyle{\sup_{x=\displaystyle{\sum_{i=1}^{n}}u_{i}\gamma_{i}v_{i}}}[\min_{1
\leq i \leq n } [\min [ \mu_{1}(u_{i}), \mu_{2}(v_{i})]] ] =
(\mu_{1} \circ \mu_{2})(x) $.
\\Thus $ \mu_{1} \circ \mu_{2}\subseteq \mu_{3} $.

\begin{Proposition} Let $ \mu_{1}, \mu_{2}, \mu_{3} \in
FLI(S)[ FRI(S), FI(S)] $. Then

(i) $ (\mu_{1} \circ \mu_{2}) \circ \mu_{3} = \mu_{1} \circ (
\mu_{2} \circ \mu_{3})$.

(ii) $ \mu_{1} \subseteq \mu_{2} $ implies that $ \mu_{1} \circ
\mu_{3} \subseteq \mu_{2} \circ \mu_{3} $.

(iii) $ \mu_{1} \circ \mu_{2} = \mu_{2} \circ \mu_{1} $, if S is
commutative $\Gamma-$semiring.

(iv) $ \textbf{1} \circ \mu_{1}= \mu_{1} $ where $ \textbf{1} \in
FLI(S) $ is defined by $ \textbf{1}(x)=1 $ for all $ x \in S $
 [resp. $ \mu_{1} \circ \textbf{1} =\mu_{1}, \textbf{1} \circ
\mu_{1} = \mu_{1} \circ \textbf{1}=
\mu_{1}]$.\label{prop:3.10}\end{Proposition}

\textsl{Proof.} Proof of (i) follows from the definition.
\\(ii) Let $ \mu_{1}\subseteq \mu_{2} $.  Now $ (\mu_{1} \circ
\mu_{3})(x)\\~~~~~= \displaystyle{\sup_{ x=
\displaystyle{\sum_{i=1}^{n}} u_{i} \gamma_{i} v_{i}}}[\min_{1 \leq
i \leq n}[\min[\mu_{1}(u_{i}), \mu_{3}(v_{i})]] : u_{i}, v_{i} \in
S, \gamma_{i} \in \Gamma, n \in Z^{+}]\\~~~~~\leq \sup[\min_{1 \leq
i \leq n}[\min[\mu_{2}(u_{i}), \mu_{3}(v_{i})]]]= (\mu_{2} \circ
\mu_{3})(x) $.
\\Thus $ \mu_{1} \circ \mu_{3} \subseteq \mu_{2} \circ \mu_{3} $.
\\(iii) $ (\mu_{1} \circ \mu_{2})(x)\\= \displaystyle{\sup_{ x=
\displaystyle{\sum_{i=1}^{n}} u_{i} \gamma_{i} v_{i}}} [\min_{1 \leq
i \leq n}[\min[ \mu_{1}(u_{i}), \mu_{2}(v_{i})]] : u_{i}, v_{i} \in
S, \gamma_{i} \in \Gamma,
 n \in
Z^{+}]\\= \displaystyle{\sup_{x=\displaystyle{\sum_{i=1}^{n}} v_{i}
\gamma_{i} u_{i}}} [\min_{1 \leq i \leq n}[\min[
 \mu_{2}(v_{i}),\mu_{1}(u_{i})]]  ]$ if S is commutative $\Gamma-$semiring\\= $ ( \mu_{2} \circ
\mu_{1})(x)$. \\Hence $ \mu_{1} \circ \mu_{2} = \mu_{2} \circ
\mu_{1} $.

(iv) As S is with left unity $ \displaystyle{\sum_{i}}[e_{i},
\delta_{i}] \in L $ which is defined by \\$ \displaystyle{\sum_{i}}
e_{i}\delta_{i}x=x$(\textit{cf. Definition}
5.1\textnormal{\cite{re:Dutta}}) for every $ x \in S $ we have, \\$
(\textbf{1} \circ
\mu_{1})(x)=\displaystyle{\sup_{x=\displaystyle{\sum_{i=1}^{n}}
u_{i} \gamma_{i} v_{i}}}[\min_{1 \leq i \leq
n}[\min[\textbf{1}(u_{i}), \mu_{1}(v_{i})]] : u_{i},v_{i} \in S,
\gamma_{i} \in \Gamma, n \in Z^{+}]
\\~~~=\sup[\min_{1 \leq i \leq n}[\min[1,
\mu_{1}(v_{i})]]]=\sup[\min_{1 \leq i \leq n} \mu_{1}[(v_{i})]]\leq
\sup[\min_{1 \leq i \leq n}
[\mu_{1}(u_{i}\gamma_{i}v_{i})]]\\~~~\leq
\mu_{1}(\displaystyle{\sum_{i=1}^{n}} u_{i}
\gamma_{i}v_{i})=\mu_{1}(x)$.
\\Therefore $(\textbf{1} \circ \mu_{1}) \subseteq \mu_{1}$.
Again $(\textbf{1} \circ \mu_{1})(x)\\~~~=
\displaystyle{\sup_{x=\displaystyle{\sum_{i=1}^{n}} u_{i}\gamma_{i}
v_{i}}}[\min_{1 \leq i \leq n}[\min[\textbf{1}(u_{i}),
\mu_{1}(v_{i})]]: u_{i},v_{i} \in S; \gamma_{i} \in \Gamma, n \in
Z^{+}]\\~~~\geq \min_{1 \leq i \leq n}[\min[\textbf{1}(e_{i}),
\mu_{1}(x)]]$ [Since $ \displaystyle{\sum_{j}} e_{i} \delta_{i} x =
x]\\~~~=\mu_{1}(x)$
\\So $ \mu_{1} \subseteq \textbf{1} \circ \mu_{1}$ and hence $
\textbf{1} \circ \mu_{1} = \mu_{1}$.
\\\\The following result shows that `.' distributive over `$\oplus$'
from both sides.
\begin{Proposition} Let $ \mu_{1}, \mu_{2}, \mu_{3} \in
FLI(S) [ FRI(S), FI(S)]$. Then \\(i) $ \mu_{1} \circ (\mu_{2} \oplus
\mu_{3})= \mu_{1} \circ \mu_{2} \oplus \mu_{1} \circ \mu_{3},$
and\\(ii) $(\mu_{2} \oplus \mu_{3}) \circ \mu_{1}= \mu_{2} \circ
\mu_{1} \oplus \mu_{3} \circ
\mu_{1}$.\label{prop:3.11}\end{Proposition}

\textsl{Proof.} Since $ \mu_{2} \subseteq \mu_{2} \oplus \mu_{3}$
therefore $ \mu_{1} \circ \mu_{2} \subseteq \mu_{1} \circ (\mu_{2}
\oplus \mu_{3})$.
\\Similarly $ \mu_{1} \circ \mu_{3} \subseteq \mu_{1} \circ (\mu_{2}
\oplus \mu_{3})$.
\\Thus $  (\mu_{1} \circ \mu_{2}) \oplus(\mu_{1} \circ \mu_{3}) \subseteq (\mu_{1} \circ (\mu_{2}
\oplus \mu_{3})) \oplus (\mu_{1} \circ (\mu_{2} \oplus
\mu_{3}))\\~~~~~~~~~~~~~~~~~~~~~~~~~~~=(\mu_{1} \circ (\mu_{2}
\oplus \mu_{3}))$.
\\Now let $ x \in S $ be arbitrary.Then
\\$[\mu_{1} \circ (\mu_{2} \oplus \mu_{3})](x)\\=\displaystyle{\sup_{x=\displaystyle{\sum_{i=1}^{n}}
u_{i} \gamma_{i} v_{i}}}[\min_{1 \leq i \leq n}[\min[\mu_{1}(u_{i}),
(\mu_{2} \oplus \mu_{3})(v_{i})]]: u_{i},v_{i} \in S, \gamma_{i} \in
\Gamma, n \in Z^{+}]\\=\sup[\displaystyle{\min_{1 \leq i \leq
n}}[\min[\mu_{1}(u_{i}),\displaystyle{\sup_{v_{i}=r_{i}+s_{i}}}[\min[\mu_{2}(r_{i}),\mu_{3}(s_{i})]]]]]
\\=\displaystyle{\sup_{x=\displaystyle{\sum_{i=1}^{n}} (u_{i} \gamma_{i} r_{i}+ u_{i}
\gamma_{i} s_{i})}}[\min_{1 \leq i \leq
n}[\min[\mu_{1}(u_{i}),\mu_{2}(r_{i}),\mu_{3}(s_{i})]]
\\\leq
\displaystyle{\sup_{x=\displaystyle{\sum_{j=1}^{n}} p_{j} \delta_{j}
q_{j} + \displaystyle{\sum_{k=1}^{m}} p_{K}^{'} \delta_{k}^{'}
q_{k}^{'}}}[\min[\min[\min_{1 \leq j
\leq}[\mu_{1}(p_{j}),\mu_{2}(q_{j})]], \min[\min_{1 \leq k \leq
m}[\mu_{1}(p_{k}^{'}),\mu_{3}(q_{k}^{'})]]]]\\=\sup[\min[(\mu_{1}
\circ \mu_{2})(u), (\mu_{1} \circ \mu_{3})(v)]:
u=\displaystyle{\sum_{j=1}^{n}} p_{j} \delta_{j} q_{j}$ and $
v=\displaystyle{\sum_{k=1}^{m}} p{k}^{'}
\delta_{k}^{'}q_{k}^{'}]\\=((\mu_{1}\circ \mu_{2}) \oplus
(\mu_{1}\circ \mu_{3}))(x)$.
\\Thus $ \mu_{1} \circ(\mu_{2}\oplus \mu_{3}) \subseteq (\mu_{1}\circ \mu_{2})
\oplus (\mu_{1}\circ \mu_{3})$.
\\Hence we conclude that $ \mu_{1} \circ(\mu_{2}\oplus \mu_{3}) = (\mu_{1}\circ \mu_{2})
\oplus (\mu_{1}\circ \mu_{3})$.
\\Proof of (ii) follows similarly.

\begin{Theorem} Let S be a $ \Gamma-$semiring. Then FLI(S) and
FRI(S)both are zero-sum free hemiring having infinite element 1
under the operations of sum and composition of fuzzy left ideals and
fuzzy right ideals respectively.\label{th:3.12}\end{Theorem}

\textsl{Proof.} It is easy to see that $\theta \in FLI(S) $. Now by
using Propositions \ref{prop:3.2}, \ref{prop:3.3}, \ref{prop:3.4},
\ref{prop:3.10}, \ref{prop:3.11} for any $ \mu_{1}, \mu_{2}, \mu_{3}
\in FLI(S) $, we easily obtain
\\(i) $ \mu_{1} \oplus \mu_{2} \in FLI(S) $,
\\(ii)$ \mu_{1} \circ \mu_{2} \in FLI(S) $,
\\(iii) $ \mu_{1} \oplus \mu_{2}= \mu_{2} \oplus \mu_{1}$,
\\(iv) $\theta \oplus \mu_{1}=\mu_{1}$,
\\(v) $ \mu_{1} \oplus (\mu_{2} \oplus \mu_{3}) = (\mu_{1} \oplus
\mu_{2}) \oplus \mu_{3} $,
\\(vi) $ \mu_{1} \circ (\mu_{2} \circ \mu_{3}) = (\mu_{1} \circ
\mu_{2}) \circ \mu_{3} $,
\\(vii) $ \mu_{1} \circ (\mu_{2} \oplus \mu_{3}) = (\mu_{1} \circ
\mu_{2}) \oplus (\mu_{1} \circ \mu_{3}) $,
\\(viii) $(\mu_{2} \oplus \mu_{3}) \circ \mu_{1} = (\mu_{2} \circ
\mu_{1}) \oplus (\mu_{3} \circ \mu_{1})$.
\\Cosequently, FLI(S) is a hemiring under the operations of sum and
composition of fuzzy ideals of S.
\\Now by Proposition 3.3(v), $\textbf{1} \subseteq \textbf{1} \oplus
\mu $ for $ \mu \in FLI(S)$.
\\Also $ (\textbf{1} \oplus \mu)(x)= \displaystyle{\sup_{ x=y+z}
[\min[\textbf{1}(y), \mu(z)]: y,z \in S]} \leq 1 = \textbf{1}(x) $
for all $ x \in S $.
\\Therefore $ \textbf{1} \oplus \mu \subseteq \textbf{1} $ and
hence $ \textbf{1} \oplus \mu = \textbf{1} $ for all $ \mu \in
FLI(S) $.
\\Thus $ \textbf{1}$ is an infinite element of FLI(S).
Now let $ \mu_{1} \oplus \mu_{2} = \theta $ for \\$ \mu_{1}, \mu_{2}
\in FLI(S) $. Then $ \mu_{1} \subseteq \mu_{1} \oplus \mu_{2} =
\theta \subseteq \mu_{1} $. Consequently, $ \mu_{1}=\theta$.
\\Similarly it can be shown that $ \mu_{2}=\theta$.
Hence the hemiring FLI(S) is zero-sum free.
\\In analogous manner we can proof the result for FRI(S).

\textbf{Remark.} If S is a commutative $\Gamma-$semiring then FLI(S)
and FRI(S) are semirings.

\begin{Corollary} FI(S) is a zero-sum free simple semiring under
the operations of sum and composition of fuzzy
ideals.\label{cor:3.13}\end{Corollary}

\textsl{Proof.} By Proposition 3.10(iv) we have $ \textbf{1} \circ
\mu = \mu \circ \textbf{1} = \mu $ for all $ \mu \in FI(S) $.
\\Hence the result follows from the above theorem.

\begin{Lemma} Intersection of a nonempty collection of fuzzy
left ideals (resp. fuzzy right ideals, fuzzy ideals ) is a fuzzy
left ideal ( resp. fuzzy right ideal, fuzzy ideal) of
S.\label{lemma:3.14}\end{Lemma}

\textsl{Proof.} Let $ \{ \mu_{i} : i \in I \} $ be a nonempty family
of fuzzy ideals of S. Let  $ x, y \in S $.
\\Then $ (\displaystyle{\bigcap_{ i \in I }\mu_{i}) (x+y) = \inf_{i \in I }[
\mu_{i}(x+y)]} \displaystyle{\geq \inf_{i \in I } [\min[ \mu_{i}(x),
\mu_{i}(y)]]}\\~~~~~~~~~~~~~~~~~~~=\displaystyle{\min[ \inf_{i \in I
} [ \mu_{i}(x)],\inf_{i \in I } [\mu_{i}(y)]]}=\displaystyle{\min[
(\bigcap_{i \in I }\mu_{i})(x), (\bigcap_{i \in I } \mu_{i})(y)]}.$
\\Again $\displaystyle{(\bigcap_{i \in I }\mu_{i})(x\gamma y)= \inf_{i \in I } [
\mu_{i}(x\gamma y)]}\displaystyle{\geq \inf_{i \in I
}[\mu_{i}(y)]=(\bigcap_{i \in I } \mu_{i})(y)}$.
\\Thus $\displaystyle{ \bigcap_{i \in I }\mu_{i}}$ is a fuzzy left ideal of S.
\\Similarly we can prove the other statements.

\begin{Theorem} Let $ \mu_{1}$ and $\mu_{2}$ be two fuzzy
left ideals (fuzzy right ideals, fuzzy ideals) of a
$\Gamma-$semiring S. Then $ \mu_{1} \oplus \mu_{2}$ is the unique
minimal element of the family of all fuzzy left ideals (resp. fuzzy
right ideals, fuzzy ideals) of S containing $\mu_{1}$ and $ \mu_{2}$
and $\mu_{1} \cap \mu_{2}$ is the unique maximal element of the
family of all fuzzy left ideals (resp. fuzzy right ideals, fuzzy
ideals) of S contained in $\mu_{1}$ and $
\mu_{2}$.\label{th:3.15}\end{Theorem}

\textsl{Proof.} Let$ \mu_{1}, \mu_{2} \in FLI(S) $. Then $ \mu_{1},
\mu_{2} \subseteq \mu_{1} \oplus \mu_{2} $ [cf. Proposition 3.3(v)].
Suppose $ \mu_{1} \subseteq \psi $ and $\mu_{2} \subseteq \psi $
where $ \psi \in FLI(S) $. Now for any $ x \in S $,
\\ $(\mu_{1} \oplus \mu_{2})(x)= \displaystyle{\sup_{x=y+z}[\min[\mu_{1}(y), \mu_{2}(z)] :
 y, z \in S]}\leq \sup[\min[\psi(y),
\psi(z)]]\\~~~~~~~~~~~~~~~~~\leq \sup \psi(y+z) = \psi(x) $
\\Thus $ \mu_{1} \oplus \mu_{2} \subseteq \psi $.
Again $ \mu_{1} \cap \mu_{2} \subseteq \mu_{1}, \mu_{2} $.
\\Let us suppose that $ \phi \in FLI(S) $ be such that $ \phi
\subseteq \mu_{1} $ and $ \phi \subseteq \mu_{2} $. Then for any $ x
\in S $,\\ $ (\mu_{1} \cap \mu_{2})(x)=\min[\mu_{1}(x), \mu_{2}(x)]
\geq \min [\phi(x), \phi(x)] =\phi(x) $.
\\Thus $ \phi \subseteq \mu_{1} \cap \mu_{2} $.
Uniqueness of $ \mu_{1} \oplus \mu_{2}$ and $ \mu_{1} \cap \mu_{2}$
with the stated properties are obvious.
\\Proofs of other cases follow similarly.

\begin{Theorem} FLI(S) [resp. FRI(S), FI(S)] is a complete
lattice.\label{th:3.16}\end{Theorem}

\textsl{Proof.} We define a relation `$\leq$' on FLI(S) as follows:
$ \mu_{1} \leq \mu_{2} $ if and only if $ \mu_{1}(x) \leq \mu_{2}(x)
$ for all $ x \in S $. Then FLI(S) is a poset with respect to
`$\leq$'.By Theorem \ref{th:3.15}, every pair of elements of FLI(S)
has lub and glb in FLI(S). Thus FLI(S) is a lattice. Now $\textbf{1}
\in FLI(S)$ and $\mu \leq \textbf{1}$ for all $\mu \in FLI(S)$. So
$\textbf{1}$ is the greatest element of FLI(S). Let $\{ \mu_{i}: i
\in I \}$ be a non empty family of fuzzy left ideals of S. Then by
Lemma \ref{lemma:3.14}, it follows that $\displaystyle{\bigcap_{i
\in I }} \mu_{i} \in FLI(S) $. Also it is the glb of $\{ \mu_{i}: i
\in I \}$. Hence FLI(S) is a complete lattice.
\\Proofs of other cases follow similarly.

\begin{Proposition} If S is a $\Gamma-$semiring then the
lattice $(FLI(S), \oplus, \cap)~~ $
\\$[(FRI(S), \oplus, \cap), (FI(S),
\oplus, \cap)] $ is modular if each of its member is a fuzzy left
k-ideal [resp. fuzzy right k-ideal, fuzzy
k-ideal].\label{prop:3.17}\end{Proposition}

\textsl{Proof.} Let us assume that every member of FLI(S) is a fuzzy
left k-ideal and $\mu_{1}, \mu_{2}, \mu_{3} \in FLI(S)$ such that
$\mu_{2} \cap \mu_{1}=\mu_{2} \cap \mu_{3}, \mu_{2} \oplus
\mu_{1}=\mu_{2} \oplus \mu_{3}$ and $\mu_{1} \subseteq \mu_{3}$.
Then for any $x \in S$,
\\ $ \mu_{1}(x)= (\mu_{1} \oplus \mu_{1})(x)=\displaystyle{\sup_{x=u+v}}[\min[\mu_{1}(u),
\mu_{1}(v)]:  u,v \in S]\\~~~\geq \sup[\min [\mu_{1}(u), (\mu_{2}
\cap \mu_{1})(v)]]=\sup[\min[\mu_{1}(u), (\mu_{2} \cap
\mu_{3})(v)]]\\~~~=\sup[\min[\mu_{1}(u), \min[\mu_{2}(v),
\mu_{3}(v)]]]=\sup[\min[\min[\mu_{1}(u), \mu_{2}(v)],
\mu_{3}(v)]]\\~~~\geq \sup[\min[\min[\mu_{1}(u), \mu_{2}(v)], \min
[\mu_{3}(u+v), \mu_{3}(u)]]] $ [Since $ \mu_{3}$ is a left
k-ideal].$\\~~~\geq \sup[\min[\min[\mu_{1}(u), \mu_{2}(v)],
\min[\mu_{3}(u+v), \mu_{1}(u)]]]\\~~~=\sup[\min[\min[\mu_{1}(u),
\mu_{2}(v)], \mu_{3}(u+v)]]\\~~~=\min[\sup[\min[\mu_{1}(u),
\mu_{2}(v)]], \sup[ \mu_{3}(u+v)]]=\min[ (\mu_{1} \oplus
\mu_{2})(x), \mu_{3}(x)]\\~~~= \min[ (\mu_{1} \oplus \mu_{3})(x),
\mu_{3}(x)]= \mu_{3}(x)$ [Since $ \mu_{3} \subseteq \mu_{1} \oplus
\mu_{3}]$.
\\Thus $ \mu_{3} \subseteq \mu_{1}$ and hence $ \mu_{1} = \mu_{3}$.
Hence $ ( FLI(S), \oplus, \cap )$ is modular.

\end{document}